\documentclass[12pt]{amsart}
\usepackage{amsmath}
\usepackage{amssymb}
\usepackage{amsfonts}
\usepackage{amsthm}
\usepackage{amscd}
\usepackage{bm}
\usepackage[dvips]{graphicx}
\newtheorem{thrm}{Theorem}
\newtheorem{coro}{Corollary}
\newtheorem{prop}{Proposition}

\newtheorem{thm}{Theorem}[section]
\newtheorem{lem}[thm]{Lemma}
\newtheorem{pro}[thm]{Proposition}
\newtheorem{cor}[thm]{Corollary}
\newtheorem{rem}[thm]{Remark}

\newcommand{\C}{\mathbf{C}}

\newcommand{\Z}{\mathbf{Z}}
\newcommand{\Q}{\mathbf{Q}}
\newcommand{\N}{\mathbf{N}}

\begin{document}
\title[On the rings of Fricke characters]
{On the rings of Fricke characters \\ of free abelian groups}
\address[Eri Hatakenaka]{Tokyo University of Agriculture and Technology,
        2-24-16, Naka-cho, Koganei-shi, Tokyo, 184-8588, Japan}
\email{hataken@cc.tuat.ac.jp}
\address[Takao Satoh]{Department of Mathematics, Faculty of Science Division II, Tokyo University of Science,
         1-3 Kagurazaka, Shinjuku, Tokyo, 162-8601, Japan}
\email{takao@rs.tus.ac.jp}
\maketitle

\begin{center}
{\sc Eri Hatakenaka}\footnote{e-address: hataken@cc.tuat.ac.jp} and {\sc Takao Satoh}\footnote{e-address: takao@rs.tus.ac.jp} \\
\end{center}

\begin{abstract}
In this paper, we determine the structure of the rings $\mathfrak{X}_{\Q}(H)$ of Fricke characters over $\Q$ of free abelian groups
$H$ of rank $n \geq 1$.
In particular, we consider the ideal $J$ in $\mathfrak{X}_{\Q}(H)$, generated by $\mathrm{tr}\,x -2$ for any $x \in H$,
and give a $\Q$-basis of each of the graded quotients $\mathrm{gr}^k(J):=J^k/J^{k+1}$ for $k \geq 1$.
Then we introduce a weight for each element of $\mathfrak{X}_{\Q}(H)$.
By using the concept of this weight, we show that $\mathfrak{X}_{\Q}(H)$ is an integral domain.
\end{abstract}
\section{Introduction}\label{S-Int}

Let $G$ be a group generated by $x_1, \ldots, x_n$. We denote by
$R(G)$ the set $\mathrm{Hom}(G, \mathrm{SL}(2,\C))$ of all $\mathrm{SL}(2,\C)$-representations of $G$. Let
$\mathcal{F}(R(G),\C)$
be the set $\{ \chi : R(G) \rightarrow \C \}$ of all complex-valued functions on $R(G)$.
Then we can regard $\mathcal{F}(R(G),\C)$ to be a $\C$-algebra in a natural way. (See Subsection {\rmfamily \ref{Ss-Des}} for details.)
For an element $x \in G$, $\mathrm{tr}\,x \in \mathcal{F}(R(G),\C)$ is defined to satisfy the equation
\[ (\mathrm{tr}\,x)(\rho) = \mathrm{tr}\, \rho(x) \]
for any $\rho \in R(G)$. Here \lq\lq $\mathrm{tr}$" in the right hand side means the trace of $2 \times 2$ matrix $\rho(x) \in \mathrm{SL}(2,\C)$.
The element $\mathrm{tr}\,x$ in
$\mathcal{F}(R(G),\C)$ is called the Fricke character of $x$.
Let $\mathfrak{X}_{\Q}(G)$ be the $\Q$-vector subspace in $\mathcal{F}(R(G),\C)$, which is generated by all $\mathrm{tr}\,x$ for $x \in G$. Then 
$\mathfrak{X}_{\Q}(G)$ has a ring structure with the multiplication of $\mathcal{F}(R(G),\C)$.
We call $\mathfrak{X}_{\Q}(G)$ the ring of Fricke characters of $G$ over $\Q$.

\vspace{0.5em}

Classically, Fricke characters were begun to study by Fricke for a free group $F_n$ generated by $x_1, \ldots, x_n$
in connection with certain problems
in the theory of Riemann surfaces. (See \cite{Fri}.) In the 1970s, Horowitz investigated algebraic properties of
$\mathfrak{X}(G)$ by using the combinatorial group theory in \cite{Ho1} and \cite{Ho2}.
In particular, he showed in \cite{Ho1} that for any $x \in G$, the Fricke character $\mathrm{tr}\,x$
can be written as a rational polynomial in $n+ \binom{n}{2}+ \binom{n}{3}$ characters
$\mathrm{tr}\, x_{i_1} x_{i_2} \cdots x_{i_l}$ for $1 \leq l \leq 3$ and $1 \leq i_1 < i_2 < \cdots < i_l \leq n$. 
Hence when we put $\mathfrak{P}$ to be a polynomial ring
\[ \Q[t_{i_1 \cdots i_l} \,|\, 1 \leq l \leq 3, \,\, 1 \leq i_1 < i_2 < \cdots < i_l \leq n], \]
then there exists a surjective homomorphism $\pi_G : \mathfrak{P} \rightarrow \mathfrak{X}_{\Q}(G)$ defined by
$t_{i_1 \cdots i_l} \mapsto \mathrm{tr}\, x_{i_1} x_{i_2} \cdots x_{i_l}$.
The study of the ring structure of $\mathfrak{X}_{\Q}(G)$ is inextricably associated with to that of the ideal
$\mathrm{Ker}(\pi_G)$. In general, however, it is quite a difficult problem to find a generating set of $\mathrm{Ker}(\pi_G)$ even
in the case that $G$ is a free group. (See also Subsection {\rmfamily \ref{Ss-Des}}.)

\vspace{0.5em}

In the case that $G$ is an {\bf abelian} group,  it is easy to see that any Fricke character of $G$ is written as a rational polynomial
in $n+ \binom{n}{2}$ characters
$\mathrm{tr}\, x_{i}$ for $1 \leq i \leq n$ and $\mathrm{tr}\, x_{i_1}x_{i_2}$ for $1 \leq i_1 < i_2 \leq n$. 
Hence when we put $\mathcal{P}$ to be a polynomial ring
\[ \Q[t_{i}, \, t_{i_1 i_2} \,|\, 1 \leq i \leq n, \,\, 1 \leq i_1 < i_2 \leq n], \]
then there exists the surjective homomorphism $\bar{\pi}_G : \mathcal{P} \rightarrow \mathfrak{X}_{\Q}(G)$
defined by $t_i \mapsto \mathrm{tr}\,x_i$ and
$t_{i_1 i_2} \mapsto \mathrm{tr}\,x_{i_1}x_{i_2}$.
We denote by $I$ the kernel of $\bar{\pi}_G$.

\vspace{0.5em}

In the present paper, we consider the case where $G$ is a free abelian group $H:=H_1(F_n,\Z)$ of rank $n$.
First, we introduce a descending filtration in $\mathfrak{X}_{\Q}(H)$.
Set $t_{i_1 \cdots i_l}' := t_{i_1 \cdots i_l} -2 \in \mathcal{P}$. We also denote by $t_{i_1 \cdots i_l}'$ its coset class in
$\mathcal{P}/I \cong \mathfrak{X}_{\Q}(H)$ by abuse of language.
Consider the ideal $J_0$ in $\mathcal{P}$ generated by all $t_{i_1 \cdots i_l}'$s. That is,
\[ J_0 = (t_i', \, t_{i_1 i_2}' \,|\, 1 \leq i \leq n, \,\, 1 \leq i_1 < i_2 \leq n) \subset \mathcal{P}. \]
Set $J:=\bar{\pi}_H(J_0)$. Then, we have a descending filtration
\[ J \supset J^2 \supset J^3 \supset \cdots \]
of ideals in $\mathcal{P}/I$. (See Subsection {\rmfamily \ref{Ss-Des}} for details.)
Set $\mathrm{gr}^k(J) := J^k/J^{k+1}$ for $k \geq 1$. Each $\mathrm{gr}^k(J)$ is a $\Q$-vector space of finite dimension
on which $\mathrm{Aut}\,H$ naturally acts.
In general, to determine the structures of the graded quotients $\mathrm{gr}^k(J)$ plays an important roles on various studies of the ring
$\mathcal{P}/I$.
The first purpose of this paper is to give a basis of $\mathrm{gr}^k(J)$ for any $k \geq 1$. More precisely, we show
the following theorem.
\begin{thrm}($=$ Theorems {\rmfamily \ref{T-marine}} and {\rmfamily \ref{T-sunshine}}.)\label{T-I-1}
\begin{enumerate}
\item For each $k \geq 1$ and $0 \leq l \leq k$, set
\[\begin{split}
   T_l := \{ t_{p_1q_1}' & \cdots t_{p_lq_l}' t_{i_{l+1}}' \cdots t_{i_k}' \in J_0 \\
   & \,|\, 1 \leq p_1 < q_1 < \cdots < p_l < q_l \leq n, \,\, 1 \leq i_{l+1} \leq \cdots \leq i_k \leq n \}.
  \end{split}\]
Then
\[ S_k := \bigcup_{l=0}^k \bar{\pi}_H(T_l) \]
forms a basis of $\mathrm{gr}^k(J)$. \\
\item $\bigcap_{k \geq 1} J^k = \{0 \}$.
\end{enumerate}
\end{thrm}

\vspace{0.5em}

Needless to say, by an argument similar to the above, we can define the ideal $J_F$ in the ring $\mathfrak{X}_{\Q}(F_n)$
of Fricke characters of the free group $F_n$ generated by all $t_{i_1 \cdots i_l}' := t_{i_1 \cdots i_l} -2$
for $1 \leq l \leq 3$ and $1 \leq i_1 < i_2 < \cdots < i_l \leq n$. Then we have a descending
filtration $J_{F_n} \supset J_{F_n}^2 \supset \cdots$ in $\mathfrak{X}_{\Q}(F_n)$, and
the graded quotients $\mathrm{gr}^k(J_{F_n}):=J_{F_n}^k/J_{F_n}^{k+1}$ for each $k \geq 1$.
Such graded quotients was originally studied by Magnus \cite{Mg2} to investigate the behavior of the action of $\mathrm{Aut}\,F_3$
on $\mathrm{gr}^k(J_{F_n})$. 
In \cite{HS1}, we gave bases of the graded quotients $\mathrm{gr}^k(J_{F_n}):=J_{F_n}^k/J_{F_n}^{k+1}$ for $k=1$ and $2$.
In general, however, it seems that no basis of $\mathrm{gr}^k(J_{F_n})$ is obtained for $k \geq 4$.
From Theorem {\rmfamily \ref{T-I-1}}, we can give a lower bound on a dimension of $\mathrm{gr}^k(J_{F_n})$ since
the natural projection $F_n \rightarrow H$ induces the surjective homomorphism $\mathrm{gr}^k(J_{F_n}) \rightarrow \mathrm{gr}^k(J)$.
More precisely, we have
\begin{coro}
For any $n \geq 2$,
\[\begin{split}
   \mathrm{dim}_{\Q}(\mathrm{gr}^k(J_{F_n})) & \geq \mathrm{dim}_{\Q}(\mathrm{gr}^k(J)) \\
    & = \sum_{l=0}^k \binom{n}{2l} \binom{n + k-l -1}{k-l}.
  \end{split}\]
\end{coro}

For any $k \geq 1$, let $\mathcal{E}_{F_n}(k)$ be the subgroup of $\mathrm{Aut}\,F_n$, consisting of automorphisms
which act on $J_{F_n}/J_{F_n}^{k+1}$ trivially.
Then we have a descending filtration $\{ \mathcal{E}_{F_n}(k) \}_{k \geq 1}$ of $\mathrm{Aut}\,F_n$.
This filtration is a Fricke character analogue of the Andreadakis-Johnson filtration of $\mathrm{Aut}\,F_n$.
The Andreadakis-Johnson filtration was originally introduced by Andreadakis \cite{And} in 1960's.
The name \lq\lq Johnson" comes from Dennis Johnson who studied this type of filtration for the
mapping class group of a surface in 1980's. It is called the Johnson filtration of the
mapping class group.
The Johnson homomorphisms are originally introduced by Johnson
in order to investigate the graded quotients of the Johnson filtration
in a series of his pioneering works \cite{Jo1}, \cite{Jo2}, \cite{Jo3} and \cite{Jo4}.
In \cite{Mo1}, Morita began to study the Johnson homomorphisms of the mapping class groups and $\mathrm{Aut}\,F_n$
systematically.
Today, together with the theory of the Johnson
homomorphisms, the Andreadskis-Johnson filtration is one of powerful tools to study the group structure of the automorphism group of a group.
(For basic material for the Andreadakis-Johnson filtration and the Johnson homomorphisms of $\mathrm{Aut}\,F_n$, see \cite{S06} or \cite{S15} for example.)
In our previous paper \cite{HS1}, we also studied the graded quotients $\mathrm{gr}^k(\mathcal{E}_{F_n}):=\mathcal{E}_{F_n}(k)/\mathcal{E}_{F_n}(k+1)$,
by introducing and using a Johnson homomorphism like homomorphism $\eta_k$ for each $k \geq 1$.
The target of the homomorphism $\eta_k$ is $\mathrm{Hom}_{\Q}(\mathrm{gr}^1(J_{F_n}), \mathrm{gr}^{k+1}(J_{F_n}))$.
In general, to determine the image of $\eta_k$ is quite a difficult problem. Since $\eta_k$ vanishes through
the homomorphism
\[ \mathrm{Hom}_{\Q}(\mathrm{gr}^1(J_{F_n}), \mathrm{gr}^{k+1}(J_{F_n})) \rightarrow \mathrm{Hom}_{\Q}(\mathrm{gr}^1(J_{F_n}), \mathrm{gr}^{k+1}(J)) \]
induced from the natural projection $F_n \rightarrow H$, from Theorem {\rmfamily \ref{T-I-1}}, we can obtain an information about the cokernel of $\eta_k$.
\begin{coro}
For any $n \geq 3$ and $k \geq 2$,
\[ \mathrm{dim}_{\Q}(\mathrm{Coker}(\eta_k)) \geq \Big{(} n + \binom{n}{2} + \binom{n}{3} \Big{)} \sum_{l=0}^{k+1} \binom{n}{2l} \binom{n + k-l}{k+1-l}. \]
\end{coro}

\vspace{0.5em}

Now, by using Theorem {\rmfamily \ref{T-I-1}}, we can obtain the following.
\begin{thrm}($=$ Theorem {\rmfamily \ref{T-blossom}}.)
The ideal $I$ is generated by 
{\small
\[\begin{split}
     t'_{ir} t'_{js} & - t'_{is} t'_{jr} \\
       & - \{ t'_i t'_r + t_j' t_s' - t_j' t_r' - t_i' t_s' \} \\
       & + \{ t_i't_{jr}' + t_j' t_{is}' + t_r' t_{is}' + t_s' t_{jr}' - t_j' t_{ir}' - t_i' t_{js}'
               - t_r' t_{js}' - t_s' t_{ir}' \} \\
       & + \frac{1}{2} \{ t_j' t_r' t_{is}' + t_i' t_s' t_{jr}' - t_i' t_r' t_{js}' - t_j' t_s' t_{ir}' \}
  \end{split} \]
}
for any $1 \leq i, j, r, s \leq n$. Here remark that in the above notation, $t_{ij}'$ should be read
\[ \begin{cases} t_{ji}' \hspace{1em} & \mathrm{if} \,\,\, i>j, \\
                 (t_i')^2 + 4 t_i' & \mathrm{if} \,\,\, i=j.
   \end{cases} \]
In particular, $I$ is finitely generated.
\end{thrm}

\vspace{0.5em}

Furthermore, from Theorem {\rmfamily \ref{T-I-1}}, we see that the ideal $I$ is generated by polynomials of degree greater than one.
This fact enables us to define the concept of a weight for each element of $\mathcal{P}/I \cong \mathfrak{X}_{\Q}(H)$.
Using this, we show
\begin{thrm} \label{I-T-P}($=$ Theorem {\rmfamily \ref{T-moonlight}}.)
The ring $\mathfrak{X}_{\Q}(H)$ of Fricke characters is an integral domain. That is, the ideal $I$ is a prime ideal in $\mathcal{P}$.
\end{thrm}

\vspace{0.5em}

Finally, in Section {\rmfamily \ref{S-Rem}}, we will give some remarks on the natural action of $\mathrm{Aut}\,H$ on $J/J^{k+1}$ for $k \geq 1$.
Let $\mathcal{E}_H(k)$ be the kernel of a natural homomorphism $\mathrm{Aut}\,H \rightarrow \mathrm{Aut}(J/J^{k+1})$ induced from the
action of $\mathrm{Aut}\,G$.
Let $\iota \in \mathrm{Aut}\,H$ be an automorphism of $H$ defined by
\[ x_i^{\iota} := x_i^{-1}, \hspace{1em} 1 \leq i \leq n. \]
Then we show the following.
\begin{prop}($=$ Proposition {\rmfamily \ref{P-erika}} and Corollary {\rmfamily \ref{C-yuri}}.)
For any $k \geq 1$, the group $\mathcal{E}_H(k)$ is the cyclic group of order $2$, generated by $\iota$.
\end{prop}

\tableofcontents

\section{Notation and conventions}\label{S-Not}

\vspace{0.5em}

Throughout the paper, we use the following notation and conventions. Let $F_n$ be the free group of rank $n$ with a basis $x_1, \ldots, x_n$,
and let $H$ be its abelianization$H_1(F_n,\Z)$. Then $H$ is a free abelian group of rank $n$, and the coset classes of $x_1, \ldots, x_n$
form a basis of $H$ as a free abelian group.
We also use the following notations.
\begin{itemize}
\item Let $G$ be a group. The automorphism group $\mathrm{Aut}\,G$ of $G$ acts on $G$ from the right.
      For any $\sigma \in \mathrm{Aut}\,G$ and $x \in G$, the action of $\sigma$ on $x$ is denoted by $x^{\sigma}$.
\item Let $N$ be a normal subgroup of a group $G$. For an element $g \in G$, we also denote the coset class of $g$ by $g \in G/N$
      if there is no confusion.
      Similarly, for a ring $R$, an element $f \in R$ and an ideal $I$ of $R$, we also denote by $f$ the coset class of $f$ in $R/I$ if
      there is no confusion.
\item For elements $x$ and $y$ in $G$, the commutator bracket $[x,y]$ of $x$ and $y$
      is defined to be $xyx^{-1}y^{-1}$.
\end{itemize}
For pairs $(i_1, i_2, \ldots, i_k)$ and $(j_1, j_2, \ldots, j_k)$ of natural numbers $i_r, j_s \in \N$, we denote the lexicographic order
among them by $(i_1, i_2, \ldots, i_k) \leq (j_1, j_2, \ldots, j_k)$. Namely, this means $i_1 < j_1$, $i_1 = j_1$ and $i_2 < j_2$, or and so on.

\vspace{0.5em}

\section{A filtration on the ring of Fricke characters of $H$}\label{S-Fri}

\vspace{0.5em}

In this section, we briefly review the ring of Fricke characters of a finitely generated group $G$.

\vspace{0.5em}

\subsection{The ring $\mathfrak{X}_{\Q}(G)$ of Fricke characters of a finitely generated group $G$}\label{Ss-Des}
\hspace*{\fill}\ 

\vspace{0.5em}

Let $G$ be a group generated by $x_1, \ldots, x_n$. We denote by
$R(G)$ the set $\mathrm{Hom}(G, \mathrm{SL}(2,\C))$ of all $\mathrm{SL}(2,\C)$-representations of $G$. Let
$\mathcal{F}(R(G),\C)$
be the set $\{ \chi : R(G) \rightarrow \C \}$ of all complex-valued functions on $R(G)$.
Then $\mathcal{F}(R(G),\C)$ has a $\C$-algebra structure by the operations defined by
\[\begin{split}
  (\chi+ \chi')(\rho) & := \chi(\rho) + \chi'(\rho), \\
  (\chi \chi')(\rho) & := \chi(\rho) \chi'(\rho), \\
  (\lambda \chi)(\rho) & := \lambda (\chi(\rho)) 
  \end{split}\]
for any $\chi$, $\chi' \in \mathcal{F}(R(G),\C)$, $\lambda \in \C$, and $\rho \in R(G)$.

\vspace{0.5em}

The automorphism group $\mathrm{Aut}\,G$ of $G$ naturally acts on $R(G)$ and $\mathcal{F}(R(G),\C)$ from the right by
\[ \rho^{\sigma}(x) := \rho(x^{\sigma^{-1}}), \hspace{1em} \rho \in R(G) \,\,\, \text{and} \,\,\, x \in G \]
and
\[ \chi^{\sigma}(\rho) := \chi(\rho^{\sigma^{-1}}), \hspace{1em} \chi \in \mathcal{F}(R(G),\C) \,\,\, \text{and} \,\,\, \rho \in R(G) \]
for any $\sigma \in \mathrm{Aut}\,G$.

\vspace{0.5em}

For any $x \in G$, we define an element $\mathrm{tr}\,x$ of $\mathcal{F}(R(G),\C)$ to be
\[ (\mathrm{tr}\,x)(\rho) := \mathrm{tr}\, \rho(x) \]
for any $\rho \in R(G)$. Here \lq\lq $\mathrm{tr}$" in the right hand side means the trace of $2 \times 2$ matrix $\rho(x)$. The element $\mathrm{tr}\,x$ in
$\mathcal{F}(R(G),\C)$ is called the Fricke character of $x \in G$. 
The action of an element $\sigma \in \mathrm{Aut}\,G$ on $\mathrm{tr}\,x$ is given by $\mathrm{tr}\,x^{\sigma}$.
We have the following well-known formulae:
\begin{eqnarray}
 & & \mathrm{tr}\,x^{-1} = \mathrm{tr}\,x, \label{eq-1} \\
 & & \mathrm{tr}\,xy = \mathrm{tr}\,yx, \label{eq-2} \\
 & & \mathrm{tr}\,xy + \mathrm{tr}\,xy^{-1} =(\mathrm{tr}\,x)(\mathrm{tr}\,y), \label{eq-3} \\
 & & \mathrm{tr}\,xyz + \mathrm{tr}\,yxz =(\mathrm{tr}\,x)(\mathrm{tr}\,yz) + (\mathrm{tr}\,y)(\mathrm{tr}\,xz) 
      + (\mathrm{tr}\,z)(\mathrm{tr}\,xy) - (\mathrm{tr}\,x)(\mathrm{tr}\,y)(\mathrm{tr}\,z), \label{eq-4} \\
 & & \mathrm{tr}\,[x, y] = (\mathrm{tr}\,x)^2 + (\mathrm{tr}\,y)^2 + (\mathrm{tr}\,xy)^2 - (\mathrm{tr}\,x)(\mathrm{tr}\,y)(\mathrm{tr}\,xy) -2, \label{eq-4.5} \\
 & & 2\mathrm{tr}\,xyzw =(\mathrm{tr}\,x)(\mathrm{tr}\,yzw) + (\mathrm{tr}\,y)(\mathrm{tr}\,zwx)
       + (\mathrm{tr}\,z)(\mathrm{tr}\,wxy) + (\mathrm{tr}\,w)(\mathrm{tr}\,xyz) \label{eq-5} \\ 
 & & \hspace{5em} + (\mathrm{tr}\,xy)(\mathrm{tr}\,zw) - (\mathrm{tr}\,xz)(\mathrm{tr}\,yw) + (\mathrm{tr}\,xw)(\mathrm{tr}\,yz) \nonumber \\
 & & \hspace{5em} - (\mathrm{tr}\,x)(\mathrm{tr}\,y)(\mathrm{tr}\,zw) - (\mathrm{tr}\,y)(\mathrm{tr}\,z)(\mathrm{tr}\,xw)
                  - (\mathrm{tr}\,x)(\mathrm{tr}\,w)(\mathrm{tr}\,yz) \nonumber \\
 & & \hspace{5em} - (\mathrm{tr}\,z)(\mathrm{tr}\,w)(\mathrm{tr}\,xy) + (\mathrm{tr}\,x)(\mathrm{tr}\,y)(\mathrm{tr}\,z)(\mathrm{tr}\,w) \nonumber
\end{eqnarray}
for any $x, y, z, w \in G$. The equations (\ref{eq-4}) and (\ref{eq-5}) are due to Vogt \cite{Vog}.
(For details, see Section 3.4 in \cite{MaR} for example.)

\vspace{0.5em}

Let $\mathfrak{X}_{\Q}(G)$ be the $\Q$-vector subspace of $\mathcal{F}(R(G),\C)$ generated by all $\mathrm{tr}\,x$ for $x \in G$.
The set $\mathfrak{X}_{\Q}(G)$ naturally has a ring structure from (\ref{eq-3}).
We call $\mathfrak{X}_{\Q}(G)$ the ring of Fricke characters of $G$ over $\Q$.
Let $\mathfrak{P}$ be a rational polynomial ring
\[ \Q[t_{i_1 \cdots i_l} \,|\, 1 \leq l \leq 3, \,\, 1 \leq i_1 < i_2 < \cdots < i_l \leq n]. \]
of $n + \binom{n}{2} + \binom{n}{3}$ indeterminates.
Consider a ring homomorphism $\pi=\pi_G : \mathfrak{P} \rightarrow \mathcal{F}(R(G),\C)$ defined by
\[ \pi(1) := \frac{1}{2} (\mathrm{tr}\, 1_G), \hspace{1em} \pi(t_{i_1 \cdots i_l}) := \mathrm{tr}\,x_{i_1} \cdots x_{i_l}. \]
We see $\mathrm{Im}(\pi) \subset \mathfrak{X}_{\Q}(G)$.
By a classical result due to Horowitz, we have
\begin{thm}[Horowitz, \cite{Ho1}]\label{T-Hor1}
For any group $G$ generated by $x_1, \ldots, x_n$, the homomorphism $\pi: \mathfrak{P} \rightarrow \mathfrak{X}_{\Q}(G)$ is surjective.
\end{thm}
More precisely, Horowitz obtained a generating set of the ring of Fricke characters of $G$ over $\Z$ in \cite{Ho1}.
Using this and (\ref{eq-5}), we can obtain the above theorem.
We should remark that in general, the structure of an ideal
\[ \mathrm{Ker}(\pi_G)
   = \{ f \in \mathfrak{P} \,|\, f(\mathrm{tr}\,\rho(x_{i_1} \cdots x_{i_l}))=0 \,\,\, \text{for any} \,\,\, \rho \in R(G) \} \]
is quite difficult. For example, its generating set is not obtained even in the case where $G$ is a free group $F_n$ in general.
Horowitz \cite{Ho1} showed that $\mathrm{Ker}(\pi_{F_n})=(0)$ for $n=1$ and $2$, and that $\mathrm{Ker}(\pi_{F_3})$ is a principal ideal generated by a quadratic polynomial
\[ t_{123}^2 - P_{123}(t) t_{123} + Q_{123}(t) \]
where
\[\begin{split}
  P_{abc}(t) & := t_{ab} t_c + t_{ac} t_b + t_{bc} t_a - t_a t_b t_c, \\
  Q_{abc}(t) & := t_a^2 + t_b^2 + t_c^2 + t_{ab}^2 + t_{ac}^2 + t_{bc}^2
                   - t_a t_b t_{ab} - t_a t_c t_{ac} - t_b t_c t_{bc} + t_{ab} t_{bc} t_{ac} -4.
 \end{split}\]
For $n \geq 4$, Whittemore \cite{Whi} showed that the ideal $\mathrm{Ker}(\pi_{F_n})$ is not principal. However,
very little is known for  $\mathrm{Ker}(\pi_{F_n})$ for general $n \geq 4$.

\vspace{0.5em}

\subsection{The ring $\mathfrak{X}_{\Q}(G)$ for an abelian group $G$}\label{Ss-Fil}
\hspace*{\fill}\ 

\vspace{0.5em}

Now, in the following, we always consider the case where $G$ is {\bf abelian}. In this case, we see
\begin{eqnarray}
 & & 2 \mathrm{tr}\,xyz = (\mathrm{tr}\,x)(\mathrm{tr}\,yz) + (\mathrm{tr}\,y)(\mathrm{tr}\,xz) 
      + (\mathrm{tr}\,z)(\mathrm{tr}\,xy) - (\mathrm{tr}\,x)(\mathrm{tr}\,y)(\mathrm{tr}\,z), \label{eq-6}
\end{eqnarray}
from (\ref{eq-4}). This shows that if $G$ is abelian, and is generated by $x_1, \ldots, x_n$, it turns out that $\mathfrak{X}_{\Q}(G)$ is generated by
$\mathrm{tr}\,1_{G}$,
\[ \mathrm{tr}\,x_i \hspace{0.5em} \text{for} \hspace{0.5em} 1 \leq i \leq n \hspace{1em} \text{and} \hspace{1em}
   \mathrm{tr}\,x_ix_j \hspace{0.5em} \text{for} \hspace{0.5em} 1 \leq i < j \leq n \]
by Theorem {\rmfamily \ref{T-Hor1}}. In other words, consider a polynomial ring
\[ \mathcal{P} := \Q[t_i, t_{i_1 i_2} \,|\, 1 \leq i \leq n, \,\,\, 1 \leq i_1 < i_2 < \leq n], \]
and a ring homomorphism $\bar{\pi}=\bar{\pi}_G : \mathcal{P} \rightarrow \mathfrak{X}_{\Q}(G)$ defined by
\[ \bar{\pi}(1) := \frac{1}{2} (\mathrm{tr}\, 1_G), \hspace{1em} \bar{\pi}(t_{i_1 \cdots i_l}) := \mathrm{tr}\,x_{i_1} \cdots x_{i_l}. \]
Then $\bar{\pi} : \mathcal{P} \rightarrow \mathfrak{X}_{\Q}(G)$ is surjective. Set $I:= \mathrm{Ker}(\bar{\pi})$. 
In this paper, we always identify $\mathcal{P}/I$ with $\mathfrak{X}_{\Q}(G)$ as a ring under the isomorphism induced from $\bar{\pi}$, and
also call each of them the ring of Fricke characters of $G$ over $\Q$.

\vspace{0.5em}

Next, set $t_{i_1 \cdots i_l}' := t_{i_1 \cdots i_l} -2 \in \mathcal{P}$. We also denote by $t_{i_1 \cdots i_l}'$ its coset class in $\mathcal{P}/I$
by abuse of language.
An element in $\mathfrak{X}_{\Q}(G)$ corresponding to $t_{i_1 \cdots i_l}' \in \mathcal{P}/I$ is
\[\begin{split}
   \mathrm{tr}'\,x_{i_1} \cdots x_{i_l} : & = (\mathrm{tr}\,x_{i_1} \cdots x_{i_l}) -2 \\
     & = (\mathrm{tr}\,x_{i_1} \cdots x_{i_l}) - \mathrm{tr}\, 1_G \in \mathfrak{X}_{\Q}(G).
  \end{split}\]
Consider an ideal
\[ J_0 = (t_i', \, t_{i_1 i_2}' \,|\, 1 \leq i \leq n, \,\, 1 \leq i_1 < i_2 \leq n) \subset \mathcal{P}. \]
generated by all $t_{i_1 \cdots i_l}'$'s in $\mathcal{P}$, and
\[ J := \bar{\pi}(J_0) = (t_i', \, t_{i_1 i_2}' \,|\, 1 \leq i \leq n, \,\, 1 \leq i_1 < i_2 \leq n) \subset \mathcal{P}/I. \]
Then, we have a descending filtration
\[ J \supset J^2 \supset J^3 \supset \cdots \]
of ideals of $\mathcal{P}/I$. Set
\[ \mathrm{gr}^k(J) := J^k/J^{k+1}. \]
Then each of $\mathrm{gr}^k(J)$ is a finite dimensional $\Q$-vector space.
In the present paper, for the case where $G$ is a free abelian group $H$ of rank $n$,
we determine the $\Q$-vector space structures of $\mathrm{gr}^k(J)$ for any $k \geq 1$.
More precisely, we give a basis of each of $\mathrm{gr}^k(J)$ in Section {\rmfamily \ref{S-Str}}.

\vspace{0.5em}

\subsection{Basic formulae among $\mathrm{tr}'\,x$}\label{Ss-For}
\hspace*{\fill}\ 

\vspace{0.5em}

In this subsection, we summarize basic and useful formulae among $\mathrm{tr}'\,x$ for $x \in G$. To begin with, we confirm the following.
\begin{eqnarray}
 & & \mathrm{tr}'\,x^{-1} = \mathrm{tr}'\,x, \label{eq-6} \\
 & & \mathrm{tr}'\,xy = \mathrm{tr}'\,yx, \label{eq-7} \\
 & & \mathrm{tr}'\,xy + \mathrm{tr}'\,xy^{-1} =2\mathrm{tr}'\,x + 2\mathrm{tr}'\,y + (\mathrm{tr}'\,x)(\mathrm{tr}'\,y), \label{eq-8} \\
 & & \mathrm{tr}'\,xyz + \mathrm{tr}'\,yxz = -2 \{ \mathrm{tr}'\,x + \mathrm{tr}'\,y + \mathrm{tr}'\,z \}
                                                  + 2 \{ \mathrm{tr}'\,xy + \mathrm{tr}'\,yz + \mathrm{tr}'\,xz \} \label{eq-9} \\
 & & \hspace{8em} + (\mathrm{tr}'\,x)(\mathrm{tr}'\,yz) + (\mathrm{tr}'\,y)(\mathrm{tr}'\,xz) 
      + (\mathrm{tr}'\,z)(\mathrm{tr}'\,xy), \nonumber \\
 & & \hspace{8em} -2 \{ (\mathrm{tr}'\,x)(\mathrm{tr}'\,y) + (\mathrm{tr}'\,y)(\mathrm{tr}'\,z) + (\mathrm{tr}'\,z)(\mathrm{tr}'\,x) \} \nonumber \\
 & & \hspace{8em} - (\mathrm{tr}'\,x)(\mathrm{tr}'\,y)(\mathrm{tr}'\,z), \nonumber \\
\end{eqnarray}
and
{\small
\begin{equation}\begin{split}
2\mathrm{tr}'\,xyzw & = 2(\mathrm{tr}'\,x + \mathrm{tr}'\,y + \mathrm{tr}'\,z + \mathrm{tr}'\,w) \\
   & \hspace{2em} - 2(\mathrm{tr}'\,xy + \mathrm{tr}'\,xz + \mathrm{tr}'\,xw + \mathrm{tr}'\,yz + \mathrm{tr}'\,yw + \mathrm{tr}'\,zw) \\
   & \hspace{2em} +2(\mathrm{tr}'\,xyz + \mathrm{tr}'\,xyw + \mathrm{tr}'\,xzw + \mathrm{tr}'\,yzw) \\
   & \hspace{2em} +2 \{ (\mathrm{tr}'\,x)(\mathrm{tr}'\,y) + (\mathrm{tr}'\,x)(\mathrm{tr}'\,w) + (\mathrm{tr}'\,y)(\mathrm{tr}'\,z)
                     + (\mathrm{tr}'\,z)(\mathrm{tr}'\,w) \\
   & \hspace{5em} + 2 (\mathrm{tr}'\,x)(\mathrm{tr}'\,z) + 2 (\mathrm{tr}'\,y)(\mathrm{tr}'\,w) \} \\
   & \hspace{2em} - 2 \{ (\mathrm{tr}'\,x)(\mathrm{tr}'\,yz) + (\mathrm{tr}'\,x)(\mathrm{tr}'\,zw) + (\mathrm{tr}'\,y)(\mathrm{tr}'\,xw)
                      + (\mathrm{tr}'\,y)(\mathrm{tr}'\,zw) \\
   & \hspace{5em} + (\mathrm{tr}'\,z)(\mathrm{tr}'\,xy) + (\mathrm{tr}'\,z)(\mathrm{tr}'\,xw) + (\mathrm{tr}'\,w)(\mathrm{tr}'\,xy)
                      + (\mathrm{tr}'\,w)(\mathrm{tr}'\,yz) \} \\
   & \hspace{2em} + \{ (\mathrm{tr}'\,x)(\mathrm{tr}'\,yzw) + (\mathrm{tr}'\,y)(\mathrm{tr}'\,xzw) + (\mathrm{tr}'\,z)(\mathrm{tr}'\,xyw)
                      + (\mathrm{tr}'\,w)(\mathrm{tr}'\,xyz) \} \\
   & \hspace{2em} + \{ (\mathrm{tr}'\,xy)(\mathrm{tr}'\,zw) - (\mathrm{tr}'\,xz)(\mathrm{tr}'\,yw) + (\mathrm{tr}'\,xw)(\mathrm{tr}'\,yz) \} \\
   & \hspace{2em} - \{ (\mathrm{tr}'\,x)(\mathrm{tr}'\,y)(\mathrm{tr}'\,zw) + (\mathrm{tr}'\,y)(\mathrm{tr}'\,z)(\mathrm{tr}'\,xw) 
                       + (\mathrm{tr}'\,x)(\mathrm{tr}'\,w)(\mathrm{tr}'\,yz) \\
   & \hspace{5em} + (\mathrm{tr}'\,z)(\mathrm{tr}'\,w)(\mathrm{tr}'\,xy) \} \\
   & \hspace{2em} + (\mathrm{tr}'\,x)(\mathrm{tr}'\,y)(\mathrm{tr}'\,z)(\mathrm{tr}'\,w) \\
   & \hspace{2em} + 2 \{ (\mathrm{tr}'\,x)(\mathrm{tr}'\,y)(\mathrm{tr}'\,z) + (\mathrm{tr}'\,x)(\mathrm{tr}'\,y)(\mathrm{tr}'\,w)
                  + (\mathrm{tr}'\,x)(\mathrm{tr}'\,z)(\mathrm{tr}'\,w) \\
   & \hspace{5em} + (\mathrm{tr}'\,y)(\mathrm{tr}'\,z)(\mathrm{tr}'\,w) \}.
\end{split} \label{eq-20} \end{equation}
}
These formulae hold for any group $G$ and $x, y, z, w \in G$.
For details, see section 4 in our previous paper \cite{HS1}.

\begin{lem}\label{L-candy}
For any $x \in G$ and $\alpha \in \Z$,
\[ \mathrm{tr}'\, x^{\alpha} \equiv \alpha^2 \mathrm{tr}'\, x \hspace{0.5em} \pmod{J^2}. \]
\end{lem}
\textit{Proof.}
It is obvious when $\alpha=0, 1$. It suffices to show $\alpha > 0$. We show this by the induction on $\alpha$.
Assume $\alpha \geq 2$. Then, substituting $x^{\alpha-1}$ and $x$ to $x$ and $y$ in (\ref{eq-8}) respectively, we have
\[ \mathrm{tr}'\, x^{\alpha} + \mathrm{tr}'\,x^{\alpha-2}
    =2\mathrm{tr}'\,x^{\alpha-1} + 2\mathrm{tr}'\,x + (\mathrm{tr}'\,x^{\alpha-1})(\mathrm{tr}'\,x), \]
and hence by the inductive hypothesis we obtain
\[\begin{split}
   \mathrm{tr}'\, x^{\alpha} & \equiv 2(\alpha-1)^2 \mathrm{tr}'\,x - (\alpha-2)^2 \mathrm{tr}'\,x
         + 2\mathrm{tr}'\,x, \\
   & \equiv \alpha^2 \mathrm{tr}'\, x \hspace{0.5em} \pmod{J^2}.
  \end{split}\]
This completes the proof of Lemma {\rmfamily \ref{L-candy}}. $\square$

\vspace{0.5em}

Here we consider additional relations among $\mathrm{tr}'\,x$ for abelian group $G$.
\begin{lem}\label{L-colon}
For an abelian group $G$ and any $x, y, z, w \in G$, we have
{\small
\[\begin{split}
   (\mathrm{tr}'\, xz)(\mathrm{tr}'\, yw) = & (\mathrm{tr}'\, xw)(\mathrm{tr}'\, yz) \\
       & + \{ (\mathrm{tr}'\, x)(\mathrm{tr}'\, z) + (\mathrm{tr}'\, y)(\mathrm{tr}'\, w)
               - (\mathrm{tr}'\, y)(\mathrm{tr}'\, z) - (\mathrm{tr}'\, x)(\mathrm{tr}'\, w) \} \\
       & - \{ (\mathrm{tr}'\, x)(\mathrm{tr}'\, yz) + (\mathrm{tr}'\, y)(\mathrm{tr}'\, xw)
               + (\mathrm{tr}'\, z)(\mathrm{tr}'\, xw) + (\mathrm{tr}'\, w)(\mathrm{tr}'\, yz) \\
       & \hspace{1.5em} - (\mathrm{tr}'\, y)(\mathrm{tr}'\, xz) - (\mathrm{tr}'\, x)(\mathrm{tr}'\, yw)
               - (\mathrm{tr}'\, z)(\mathrm{tr}'\, yw) - (\mathrm{tr}'\, w)(\mathrm{tr}'\, xz) \} \\
       & - \frac{1}{2} \{ (\mathrm{tr}'\, y)(\mathrm{tr}'\, z)(\mathrm{tr}'\, xw) 
            + (\mathrm{tr}'\, x)(\mathrm{tr}'\, w)(\mathrm{tr}'\, yz) - (\mathrm{tr}'\, x)(\mathrm{tr}'\, z)(\mathrm{tr}'\, yw) \\
       & \hspace{1.5em} - (\mathrm{tr}'\, y)(\mathrm{tr}'\, w)(\mathrm{tr}'\, xz) \}
  \end{split}\]
}
\end{lem}
\textit{Proof.}
In order to obtain the equation above, it suffices to calculate
\[ 2 \mathrm{tr}'\,xyzw - 2 \mathrm{tr}'\,yxzw \]
with the equation (\ref{eq-20}). The calculation is straightforward. We leave it to the reader for an exercise. $\square$

\vspace{0.5em}

As special cases of Lemma {\rmfamily \ref{L-colon}}, we see the following corollaries.

\begin{cor}\label{C-kofre}
For an abelian group $G$ and $x, y, z \in G$, we have
{\small
\[\begin{split}
   (\mathrm{tr}'\, xy)(\mathrm{tr}'\, xz) = & -3 (\mathrm{tr}'\, x)(\mathrm{tr}'\, y) -3 (\mathrm{tr}'\, x)(\mathrm{tr}'\, z)
           - (\mathrm{tr}'\,x)^2 - (\mathrm{tr}'\, y)(\mathrm{tr}'\, z) + 2 (\mathrm{tr}'\, x)(\mathrm{tr}'\, yz) \\
       &  + (\mathrm{tr}'\, x)(\mathrm{tr}'\, xy)
          + (\mathrm{tr}'\, x)(\mathrm{tr}'\, xz) + (\mathrm{tr}'\, y)(\mathrm{tr}'\, xz) + (\mathrm{tr}'\, z)(\mathrm{tr}'\, xy) \\
       & +\frac{1}{2} (\mathrm{tr}'\, x)^2(\mathrm{tr}'\, yz) - (\mathrm{tr}'\, x)^2(\mathrm{tr}'\, y)
               - (\mathrm{tr}'\, x)^2(\mathrm{tr}'\, z) \\
       & - (\mathrm{tr}'\, y)(\mathrm{tr}'\, xz) - (\mathrm{tr}'\, x)(\mathrm{tr}'\, yw)
               - (\mathrm{tr}'\, z)(\mathrm{tr}'\, yw) - (\mathrm{tr}'\, w)(\mathrm{tr}'\, xz)  \\
       & - \frac{1}{2} (\mathrm{tr}'\, x)^2(\mathrm{tr}'\, y)(\mathrm{tr}'\, z) 
            -2 (\mathrm{tr}'\, x)(\mathrm{tr}'\, y)(\mathrm{tr}'\, z) 
            + \frac{1}{2} (\mathrm{tr}'\, x)(\mathrm{tr}'\, z)(\mathrm{tr}'\, xy) \\
       & + \frac{1}{2} (\mathrm{tr}'\, x)(\mathrm{tr}'\, y)(\mathrm{tr}'\, xz)
  \end{split}\]
}
\end{cor}

\begin{cor}\label{C-sipre}
For an abelian group $G$ and $x, y \in G$, we have
\[\begin{split}
   (\mathrm{tr}'\,xy)^2 = & - (\mathrm{tr}'\,x)^2 - (\mathrm{tr}'\,y)^2 \\
    & + 2 \{ (\mathrm{tr}'\, x)(\mathrm{tr}'\, y) + (\mathrm{tr}'\, x)(\mathrm{tr}'\, xy) + (\mathrm{tr}'\, y)(\mathrm{tr}'\, xy) \} \\
    & + (\mathrm{tr}'\, x)(\mathrm{tr}'\, y)(\mathrm{tr}'\, xy)
  \end{split}\]
\end{cor}

\vspace{0.5em}

\section{The structure of a $\Q$-vector space $\mathrm{gr}^k(J)$ for $G=H$}\label{S-Str}
\hspace*{\fill}\ 

\vspace{0.5em}

The goal of this section is to give a basis of $\mathrm{gr}^k(J)$ for any $k \geq 1$ as a $\Q$-vector space.

\begin{pro}\label{P-popli}
For any $k \geq 1$, consider a polynomial
\[ f:= t_{p_1q_1}' \cdots t_{p_lq_l}'t_{i_{l+1}}' \cdots t_{i_k}' \in \mathcal{P} \]
for $0 \leq l \leq k$, $(p_1, q_1) \leq  \cdots \leq (p_l,q_l)$ and $1 \leq i_{l+1} \leq \cdots \leq i_k \leq n$.
Then, under the modulo $I$, the monomial $f$ can be written as a sum of monomials of type
\[ t_{p_1' q_1'}' \cdots t_{p_m' q_m'}'t_{j_{m+1}}' \cdots t_{j_k}' \in \mathcal{P} \]
such that
\begin{itemize}
\item $0 \leq m \leq k$,
\item $1 \leq p_1' < q_1' < \cdots < p_m' < q_m' \leq n$,
\item $1 \leq j_{m+1} \leq \cdots \leq j_k \leq n$,
\item $\{p_1', q_1', \ldots, p_m', q_m' \} \subset \{ p_1, q_1, \ldots, p_l, q_l \}$.
\end{itemize}
\end{pro}
\textit{Proof.}
We prove this proposition by the induction on $l \geq 0$. If $l=0$ or $1$, it is clear. Assume $l \geq 2$.

\vspace{0.5em}

First, assume that some elements in $\{p_1, q_1, \ldots, p_l, q_l \}$ are equal.
For simplicity, we consider three cases:
\begin{enumerate}
\item $p_1=p_2$ and $q_1 \neq q_2$,
\item $p_1=q_2$ and $q_1 \neq p_2$,
\item $p_1=p_2$ and $q_1=q_2$.
\end{enumerate}
For the parts (1) and (2),
by using Corollary {\rmfamily \ref{C-kofre}}, we see that under the modulo $I$, the monomial $f$ can be written as a sum of
monomials
\[ t_{p_1q_1}' \cdots t_{p_rq_r}'t_{i_{r+1}}' \cdots t_{i_k}' \]
such that $r < l$. Hence by applying the inductive hypothesis to each of such monomials, we obtain the required result.
We can discuss an argument similar to the above for the case (3) by using Corollary {\rmfamily \ref{C-sipre}}. 
Therefore we assume that $p_1, q_1, \ldots, p_l, q_l$ are distinct.

\vspace{0.5em}

Now, by the assumption $(p_1, q_1) \leq (p_2, q_2)$, we have $p_1<p_2$. If $q_1< p_2$, then we have inequalities:
\[\begin{split}
   p_1 < q_1 < & p_2 < p_3 < \cdots < p_l \\
               & \rotatebox{270}{$<$} \hspace{1.8em} \rotatebox{270}{$<$} \hspace{4em} \rotatebox{270}{$<$} \\
               & q_2 \hspace{1.5em} q_3 \hspace{3.7em} q_l
 \end{split}\]
Hence if we set $f' := t_{p_2q_2}' \cdots t_{p_lq_l}'t_{i_{l+1}}' \cdots t_{i_k}' \in \mathcal{P}$, then we can apply the inductive
hypothesis to $f'$, and obtain the required result.

\vspace{0.5em}

Next, if $p_2 < q_1$, we have inequalities:
\[\begin{split}
   p_1 \,\, < & \,\, p_2 \,\, < \,\, p_3  \,\, < \cdots < \,\, p_l \\
         & \,\, \rotatebox{270}{$<$} \,\, \hspace{1.8em} \,\, \rotatebox{270}{$<$} \hspace{4.7em} \rotatebox{270}{$<$} \\
        q_1&, \, q_2 \hspace{1.7em} q_3 \hspace{4.6em} q_l
 \end{split}\]
Then by using Lemma {\rmfamily \ref{L-colon}}, we can write $f$ as a sum of monomials
\[ t_{p_1p_2}' t_{q_1,q_2}' t_{p_3q_3}' \cdots t_{p_lq_l}'t_{i_{l+1}}' \cdots t_{i_k}' \]
and
\[ t_{p_1q_1}' \cdots t_{p_rq_r}'t_{i_{r+1}}' \cdots t_{i_k}' \]
such that $r < l$. If we apply the inductive hypothesis to $t_{q_1,q_2}' t_{p_3q_3}' \cdots t_{p_lq_l}'t_{i_{l+1}}' \cdots t_{i_k}'$
and $ t_{p_1q_1}' \cdots t_{p_rq_r}'t_{i_{r+1}}' \cdots t_{i_k}'$, we obtain the required result.
This completes the proof of Proposition {\rmfamily \ref{P-popli}}. $\square$

\begin{thm}\label{T-marine}
For each $k \geq 1$ and $0 \leq l \leq k$, set
\[\begin{split}
   T_l := \{ t_{p_1q_1}' \cdots t_{p_lq_l}' & t_{i_{l+1}}' \cdots t_{i_k}' \in J_0 \\
   & \,|\, 1 \leq p_1 < q_1 < \cdots < p_l < q_l \leq n, \,\, 1 \leq i_{l+1} \leq \cdots \leq i_k \leq n \}.
  \end{split}\]
Then
\[ S_k := \bigcup_{l=0}^k \bar{\pi}(T_l) \]
forms a basis of $\mathrm{gr}^k(J)$.
\end{thm}
\textit{Proof.}
By Proposition {\rmfamily \ref{P-popli}}, we see that $S_k$ generates $\mathrm{gr}^k(J)$. In order to show the linearly
independentness of the elements of $S_k$, set
\[\begin{split}
   f := \sum_{l=0}^k \sum_{\substack{p_1<q_1<\cdots<p_l<q_l \\[1pt] i_{l+1} \leq \cdots \leq i_k}} 
        a_{p_1q_1,\ldots, p_lq_l, i_{l+1}, \ldots, i_{k}} t_{p_1q_1}' \cdots t_{p_lq_l}'t_{i_{l+1}}' \cdots t_{i_k}' \in J_0^k
  \end{split}\]
for $a_{p_1q_1,\ldots, p_lq_l, i_{l+1}, \ldots, i_{k}} \in \Q$, and assume $\bar{\pi}(f) \in J^{k+1}$.

\vspace{0.5em}

Consider the interior
\[ D := \{ z \in \C \,|\, z \overline{z} < 1 \} \]
of the unit disk in $\C$. For any $s_1, \ldots, s_n \in D$, 
we define a representation $\rho : H \rightarrow \mathrm{SL}(2,\C)$ by
\[ \rho(x_{i}) := \begin{pmatrix} 1-s_i & 0 \\ 0 & (1-s_i)^{-1} \end{pmatrix} \]
for any $1 \leq i \leq n$. If we consider the power series expansion
\[ \frac{1}{1-s_i} = 1+s_i + s_i^2+ s_i^3+ \cdots \]
at the origin on $D$,
we can write each of $\mathrm{tr}'\, \rho(x_i)$ and $\mathrm{tr}'\, \rho(x_ix_j)$ as a convergent power series:
\[\begin{split}
   \mathrm{tr}'\, \rho(x_i) & = \frac{s_i^2}{1-s_i} = s_i^2 + s_i^3 + s_i^4 + \cdots,  \\
   \mathrm{tr}'\, \rho(x_i x_j) & = s_i^2 + 2 s_is_j + s_j^2 + (\text{terms of degree $\geq 3$ }). \\
  \end{split}\]

\vspace{0.5em}

Then we have
{\small
\[\begin{split}
  \bar{\pi}(f)(\rho) & = \sum_{l=0}^k \sum_{\substack{p_1<q_1<\cdots<p_l<q_l \\[1pt] i_{l+1} \leq \cdots \leq i_k}} \\
    & \hspace{2em} \Big{\{} a_{p_1q_1,\ldots, p_lq_l, i_{l+1}, \ldots, i_{k}}
     (s_{p_1}^2 + 2 s_{p_1} s_{q_1} + s_{q_1}^2) \cdots (s_{p_2}^2 + 2 s_{p_2} s_{q_2} + s_{q_2}^2) s_{i_{l+1}}^2 \cdots s_k^2 \Big{\}} \\
    & \hspace{1em} + (\text{terms of degree $\geq 2k+2$ }).
  \end{split}\]
}
Since $\bar{\pi}(f) \in J^{k+1}$, if we regard $\bar{\pi}(f)(\rho)$ as a polynomial on $s_i$s, its degree must be greater than or equal to $2k+2$.
Hence all coefficients of degree $2k$ are zero.

\vspace{0.5em}

To begin with, for any $1 \leq p_1 < q_2 < \cdots < p_k < q_k \leq n$,
by observing the coefficients of $s_{p_1}s_{q_1} \cdots s_{p_k}s_{q_k}$, we see
$a_{p_1q_1, \ldots, p_kq_k}=0$. For any $0 \leq l \leq k$, assume that
\[ a_{p_1q_1,\ldots, p_mq_m, i_{m+1}, \ldots, i_{k}}=0 \]
for any $l \leq m$.
Then for any $1 \leq p_1 < q_2 < \cdots < p_{l-1} < q_{l-1} \leq n$ and $1 \leq i_l \leq \cdots \leq i_k \leq n$,
we see
\[ a_{p_1q_1,\ldots, p_{l-1}q_{l-1}, i_{l}, \ldots, i_{k}}=0. \]
Therefore by the inductive argument, we verify that all coefficients of $f$ are equal to zero.
This shows that elements in $S_k$ are linearly independent.
This completes the proof of Theorem {\rmfamily \ref{T-marine}}. $\square$

\vspace{0.5em}
By observing the proof of the Theorem {\rmfamily \ref{T-marine}}, we have the following.

\begin{thm}\label{T-blossom}
The ideal $I$ is generated by 
{\small
\begin{equation} \begin{split}
     t'_{ir} t'_{js} & - t'_{is} t'_{jr} \\
       & - \{ t'_i t'_r + t_j' t_s' - t_j' t_r' - t_i' t_s' \} \\
       & + \{ t_i't_{jr}' + t_j' t_{is}' + t_r' t_{is}' + t_s' t_{jr}' - t_j' t_{ir}' - t_i' t_{js}'
               - t_r' t_{js}' - t_s' t_{ir}' \} \\
       & + \frac{1}{2} \{ t_j' t_r' t_{is}' + t_i' t_s' t_{jr}' - t_i' t_r' t_{js}' - t_j' t_s' t_{ir}' \}
  \end{split} \label{eq-gen} \end{equation}
}
for any $1 \leq i, j, r, s \leq n$. Here remark that in the above notation, $t_{ij}'$ should be read
\[ \begin{cases} t_{ji}' \hspace{1em} & \mathrm{if} \,\,\, i>j, \\
                 (t_i')^2 + 4 t_i' & \mathrm{if} \,\,\, i=j.
   \end{cases} \]
In particular, $I$ is finitely generated.
\end{thm}
\textit{Proof.}
Let $I'$ be an ideal of $\mathcal{P}$ generated by elements (\ref{eq-gen}). From Lemma {\rmfamily \ref{L-colon}},
we have $I' \subset I$. For any $f \in I$, observing the proof of Proposition {\rmfamily \ref{P-popli}} we see that
$f$ can be written as
\begin{equation}\begin{split}
   f \equiv \sum_{k \geq 0} \sum_{l=0}^k \sum_{\substack{p_1<q_1<\cdots<p_l<q_l \\[1pt] i_{l+1} \leq \cdots \leq i_k}} 
        a_{p_1q_1,\ldots, p_lq_l, i_{l+1}, \ldots, i_{k}} t_{p_1q_1}' \cdots t_{p_lq_l}'t_{i_{l+1}}' \cdots t_{i_k}' \in J_0^k
  \end{split}\label{eq-normal}\end{equation}
modulo $I'$. Here, in the sum of the right hand side of the equation above, $k$ runs over finitely many non-negative integers.
Hence by Theorem {\rmfamily \ref{T-marine}}, we see that all coefficients of $f$ are zero, and obtain $f \in I'$. 
This shows $I=I'$. $\square$

\begin{rem}
We remark that if $n=1$, we have $I=(0)$.
\end{rem}

\begin{thm}\label{T-sunshine}
$\bigcap_{k \geq 1} J^k = \{0 \}$.
\end{thm}
\textit{Proof.}
For any $f \in \bigcap_{k \geq 1} J^k$, we can write $f$ as (\ref{eq-normal}). Then observing the coset class of $f$ in
$\mathrm{gr}^1(J)$, we see that $a_{p_1q_1}=a_{i_1}=0$ for any $1 \leq p_1<q_1 \leq n$ and $1 \leq i_1 \leq n$.
Next, observing the coset class of $f$ in $\mathrm{gr}^2(J)$, we see that all coefficients of $f$ of degree $2$ are zero.
By repeating this argument inductively, we obtain $f=0$. This completes the proof of Theorem {\rmfamily \ref{T-sunshine}}. $\square$

\section{The primeness of the ideal $I$}\label{S-prim}

\vspace{0.5em}

In this section, we show that the ideal $I$ is a prime ideal of $\mathcal{P}$, in other words $\mathcal{P}/I$ is an integral domain.
In order to show this, we introduce a weight of an element $f \in \mathcal{P}/I$.
From Theorem {\rmfamily \ref{T-sunshine}}, for any $f \in \mathcal{P}/I \setminus \{ 0\}$,
there exists some integer $k \geq 0$ such that
$f \in J^k \setminus J^{k+1}$. Then we call $k$ the weight of $f$, and denote it by $\mathrm{wt}(f)$.

\begin{pro}\label{P-flower}
For any $f, g \in \mathcal{P}/I \setminus \{ 0\}$, $\mathrm{wt}(fg)=\mathrm{wt}(f)+\mathrm{wt}(g)$.
\end{pro}
\textit{Proof.}
It is obvious if $\mathrm{wt}(f)=0$ or $\mathrm{wt}(g)=0$. Hence we may assume $\mathrm{wt}(f), \mathrm{wt}(g) \geq 1$.
Set $k_1 := \mathrm{wt}(f)$ and $k_2:=\mathrm{wt}(g)$. Since it is clear that $\mathrm{wt}(fg) \geq k_1 + k_2$, assume
$\mathrm{wt}(fg) > k_1 + k_2$.

\vspace{0.5em}

For any integers $\alpha_1, \ldots, \alpha_n \in \Z$, consider a group homomorphism
$\rho_{(\alpha_1, \ldots, \alpha_n)} : H \rightarrow \Z = \langle \, x_1 \, \rangle$ defined by $x_i \mapsto x_1^{\alpha_i}$
for any $1 \leq i \leq n$.
Then $\rho_{(\alpha_1, \ldots, \alpha_n)}$ induces a ring homomorphism
$\overline{\rho}_{(\alpha_1, \ldots, \alpha_n)} : \mathfrak{X}_{\Q}(H) \rightarrow \mathfrak{X}_{\Q}(\Z)$ defined by
\[ \mathrm{tr}'\,x \mapsto \mathrm{tr}'\,\rho_{(\alpha_1, \ldots, \alpha_n)}(x). \]
Namely, we have
\[\begin{split}
   \overline{\rho}_{(\alpha_1, \ldots, \alpha_n)}(\mathrm{tr}'\,x_i) & = \mathrm{tr}'\,x_1^{\alpha_i}, \\
   \overline{\rho}_{(\alpha_1, \ldots, \alpha_n)}(\mathrm{tr}'\,x_ix_j) & = \mathrm{tr}'\,x_1^{\alpha_i + \alpha_j}, \\
  \end{split}\]

\vspace{0.5em}

Now, set 
\[\begin{split}
   f & := \sum_{l=0}^{k_1} \sum_{\substack{p_1<q_1<\cdots<p_l<q_l \\[1pt] i_{l+1} \leq \cdots \leq i_{k_1}}} 
        a_{p_1q_1,\ldots, p_lq_l, i_{l+1}, \ldots, i_{k_1}} t_{p_1q_1}' \cdots t_{p_lq_l}'t_{i_{l+1}}' \cdots t_{i_{k_1}}' \\
     & \hspace{1em} +  (\text{terms of degree $> k_1$ }), \\
   g & := \sum_{m=0}^{k_2} \sum_{\substack{p_1'<q_1'<\cdots<p_m'<q_m' \\[1pt] i_{m+1}' \leq \cdots \leq i_{k_2}'}} 
        a_{p_1'q_1',\ldots, p_m'q_m', i_{m+1}', \ldots, i_{k_2}'} t_{p_1'q_1'}' \cdots t_{p_m'q_m'}'t_{i_{m+1}'}' \cdots t_{i_{k_2}'}' \\
     & \hspace{1em} +  (\text{terms of degree $> k_2$ }),
  \end{split}\]
and set
$F:=\overline{\rho}_{(\alpha_1, \ldots, \alpha_n)}(f)$ and $G:=\overline{\rho}_{(\alpha_1, \ldots, \alpha_n)}(g)$.
Then $\mathrm{wt}(FG) > k_1+k_2$ in $\mathfrak{X}_{\Q}(\Z)$.
Hence the coefficient $P(\alpha_1, \ldots, \alpha_n)$ of $FG$ of degree $k_1+k_2$ is equal to zero. Here using Lemma {\rmfamily \ref{L-candy}},
we have $P(\alpha_1, \ldots, \alpha_n)=P_1P_2$ for
{\small
\[\begin{split}
   P_1 &:= \Big{(} \sum_{l=0}^{k_1} \sum_{\substack{p_1<q_1<\cdots<p_l<q_l \\[1pt] i_{l+1} \leq \cdots \leq i_{k_1}}}
     a_{p_1q_1,\ldots, p_lq_l, i_{l+1}, \ldots, i_{k_1}} (\alpha_{p_1}+\alpha_{q_1})^2 \cdots (\alpha_{p_l}+\alpha_{q_l})^2
        \alpha_{i_{l+1}}^2 \cdots \alpha_{i_{k_1}}^2 \Big{)} \\
   P_2 &:= \Big{(} \sum_{m=0}^{k_2} \sum_{\substack{p_1'<q_1'<\cdots<p_m'<q_m' \\[1pt] i_{m+1}' \leq \cdots \leq i_{k_2}'}}
     a_{p_1'q_1',\ldots, p_m'q_m', i_{m+1}', \ldots, i_{k_2}'} (\alpha_{p_1'}+\alpha_{q_1'})^2 \cdots (\alpha_{p_m'}+\alpha_{q_m'})^2
        \alpha_{i_{m+1}'}^2 \cdots \alpha_{i_{k_2}'}^2 \Big{)}.
  \end{split}\]
}

\vspace{0.5em}

Consider $P(\alpha_1, \ldots, \alpha_n)$ as a polynomial in $\Q[\alpha_1, \ldots, \alpha_n]$.
Since $P(\alpha_1, \ldots, \alpha_n)=0$ if $\alpha_1, \ldots, \alpha_n$ run over all integers, we see that
$P(\alpha_1, \ldots, \alpha_n)=0$ as a polynomial in $\Q[\alpha_1, \ldots, \alpha_n]$. Since $\Q[\alpha_1, \ldots, \alpha_n]$ is a domain,
we have $P_1=0$ or $P_2=0$ in $\Q[\alpha_1, \ldots, \alpha_n]$.

\vspace{0.5em}

Assume $P_1=0$. Then each coefficient of monomials in $\alpha_1, \ldots, \alpha_n$ in $P_1$ is equal to zero.
First, by observing the coefficient of $\alpha_{p_1} \alpha_{q_1} \cdots \alpha_{p_{k_1}} \alpha_{q_{k_1}}$, we see
\[ a_{p_1q_1, \ldots, p_{k_1}q_{k_1}} =0. \]
Furthermore by an argument similar to that in the proof of Theorem {\rmfamily \ref{T-marine}}, we obtain that
all $a_{p_1q_1,\ldots, p_lq_l, i_{l+1}, \ldots, i_{k_1}}$s are equal to zero. Hence $f \in J^{k_1+1}$.
This contradicts to $\mathrm{wt}(f)=k_1$. By the same argument, if $P_2=0$, we have $g \in J^{k_2+1}$, and a contradiction.
Thus, we conclude that $\mathrm{wt}(fg)= k_1 + k_2$. This completes the proof of Proposition {\rmfamily \ref{P-flower}}.
$\square$

\vspace{0.5em}

As a corollary, we obtain

\begin{thm}\label{T-moonlight}
The quotient ring $\mathcal{P}/I$ is an integral domain. That is, the ideal $I$ is a prime ideal in $\mathcal{P}$.
\end{thm}

\section{Some remarks}\label{S-Rem}

Finally, we give some remarks about an action of the automorphism group $\mathrm{Aut}\,H$ of $H$ on the ring $\mathfrak{X}_{\Q}(H)$.
In general, for any group $G$, the automorphism group $\mathrm{Aut}\,G$ of $G$ naturally acts on $\mathfrak{X}_{\Q}(G)$ from the right.
(See our previous paper \cite{HS1} for details.)
In particular, for any $\sigma \in \mathrm{Aut}\,G$ and $x \in G$,
the action of $\sigma \in \mathrm{Aut}\,G$ on $\mathrm{tr}'\,x \in \mathfrak{X}_{\Q}(G)$ is given by
\[ (\mathrm{tr}'\,x)^{\sigma} = \mathrm{tr}'\,x^{\sigma}. \]
Clearly, we see that the ideal $J$ generated by all $\mathrm{tr}'\,x$ for $x \in G$
is an $\mathrm{Aut}\,G$-invariant ideal of $\mathfrak{X}_{\Q}(G)$.
Hence $\mathrm{Aut}\,G$ naturally acts on $\mathrm{gr}^k(J)$ for each $k \geq 1$.

\vspace{0.5em}

Let $\mathcal{E}_G(k)$ be the kernel of a natural homomorphism $\mathrm{Aut}\,G \rightarrow \mathrm{Aut}(J/J^{k+1})$ induced from the
action of $\mathrm{Aut}\,G$. Then the groups $\mathcal{E}_G(k)$ define a descending filtration
\[ \mathcal{E}_G(1) \supset \mathcal{E}_G(2) \supset \cdots \supset \mathcal{E}_G(k) \supset \cdots \]
of $\mathrm{Aut}\,G$. In \cite{HS1}, we have showed that this filtration is central. That is
$[\mathcal{E}_G(k), \mathcal{E}_G(l)] \subset \mathcal{E}_G(k+l)$ for any $k, l \geq 1$.
This is a Fricke character analogue of the Andreadakis-Johnson filtration of $\mathrm{Aut}\,G$. (For details for the Andreadakis-Johnson filtration,
see \cite{S06} or \cite{S14}, for example.)

\vspace{0.5em}

Here, we determine $\mathcal{E}_H(1)$ and show $\mathcal{E}_H(k)=\mathcal{E}_H(1)$
for any $k \geq 1$.
Let $\iota \in \mathrm{Aut}\,H$ be an automorphism of $H$ defined by
\[ x_i^{\iota} := x_i^{-1}, \hspace{1em} 1 \leq i \leq n. \]
Then we have

\begin{pro}\label{P-erika}
$\mathcal{E}_H(1) = \langle \iota \rangle$. Namely, $\mathcal{E}_H(1)$ is the cyclic group of order $2$, generated by $\iota$.
\end{pro}
\textit{Proof.}
Clearly, we see $\mathcal{E}_H(1) \supset \langle \iota \rangle$.
For any $\sigma \in \mathcal{E}_H(1)$, set
\[ x_i^{\sigma} := x_1^{e_1(i)} x_2^{e_2(i)} \cdots x_n^{e_n(i)}, \hspace{1em} 1 \leq i \leq n. \]
It suffices to show
\begin{itemize}
\item $e_j(i) =0$ if $j \neq i$,
\item $e_1(1) = \cdots =e_n(n) = \pm 1$.
\end{itemize}

\vspace{0.5em}

For any $1 \leq i \leq n$ and any $s_i \in D$, consider a homomorphism $\rho_i : H \rightarrow \mathrm{SL}(2, \C)$
defined by
\[ {\rho_i}(x_j) := \begin{cases}
                     \begin{pmatrix} 1-s_i & 0 \\ 0 & (1-s_i)^{-1} \end{pmatrix}, \hspace{1em} & \text{if} \hspace{1em} j=i, \\
                     E_2, & \text{if} \hspace{1em} j \neq i.
                   \end{cases}\]
Then from
\[ \mathrm{tr}'\,x_i \equiv \mathrm{tr}'\,x_i^{\sigma} = \mathrm{tr}'\,x_1^{e_1(i)} x_2^{e_2(i)} \cdots x_n^{e_n(i)}
    \pmod{J^2}, \]
by substituting $\rho_i$, we see
\[ s_i^2 \equiv e_i(i)^2 s_i^2 \pmod{(s_i^3)}, \]
and hence $e_i(i)=\pm1$.
On the other hand, by substituting $\rho_j$ for $j \neq i$, we see
\[ 0 \equiv e_j(i)^2 s_j^2 \pmod{(s_j^3)}, \]
and hence $e_j(i)=0$.

\vspace{0.5em}

In order to show $e_1(1)=\cdots = e_n(n)$, assume that $e_i(i) \neq e_j(j)$ for some $i$ and $j$.
This means $(e_i(i), e_j(j))=(\pm1, \mp1)$.
For any $s_i, s_j \in D$,
consider a homomorphism $\rho_{ij} : H \rightarrow \mathrm{SL}(2,\C)$ defined by
\[ \rho_{ij}(x_r) := \begin{cases}
                         \rho_r(x_r), \hspace{1em} & \text{if} \hspace{1em} r=i, j, \\
                         E_2, \hspace{1em} & \text{if} \hspace{1em} \text{otherwise}.
                      \end{cases}\]
Then by substituting $\rho_{ij}$ to an equation
$\mathrm{tr}'x_ix_j \equiv \mathrm{tr}'\,(x_ix_j)^{\sigma} = \mathrm{tr}'\,x_i^{e_i(i)}x_j^{e_j(j)} \pmod{J^2}$,
we obtain
\[ s_i^2 +2s_is_j + s_j^2 \equiv s_i^2 -2s_is_j +s_j^2 \pmod{(s_i^3, s_i^2s_j, s_is_j^2, s_j^3)}. \]
This is a contradiction. Therefore we obtain the required result. 
This completes the proof of Proposition {\rmfamily \ref{P-erika}}. $\square$

\begin{cor}\label{C-yuri}
For any $k \geq 2$, $\mathcal{E}_H(k)=\mathcal{E}_H(1)$.
\end{cor}
\textit{Proof.}
In general, we have $\mathcal{E}_H(k) \subset \mathcal{E}_H(1)$. On the other hand, $\mathcal{E}_H(k) \supset \mathcal{E}_H(1)$
immediately follows from Proposition {\rmfamily \ref{P-erika}}. $\square$

\section{Acknowledgments}\label{S-Ack}

The second author is supported by a Grant-in-Aid for Young Scientists (B) by JSPS.

\end{document}